\documentclass[12pt]{article}
\usepackage{latexsym, amscd, amsmath, graphics, rotating, array}
\usepackage{amssymb,amsmath,multibox}
\makeindex
\title{Zeta functions of nilpotent groups \\ - singular Pfaffians}
\date{30 September 2003}
\author{Christopher Voll\thanks{Mathematical Institute, Oxford, England.}}

\newenvironment{acknowledgements}{\bigskip\noindent\textsl{Acknowledgements.}\rm}
\newenvironment{example}{\bigskip\noindent\textbf{Example. }\,\rm}
{\bigskip}

\newenvironment{proof}{\bigskip\noindent\textbf{Proof.}\rm}{\hfill$\Box$ \\ \hspace{0.1cm}}

\newenvironment{proofnodot}{\bigskip\noindent\textbf{Proof}\rm}{\hfill$\Box$ \\ \hspace{0.1cm}}

\newtheorem{lemma}{Lemma}
\newtheorem{theorem}{Theorem}

\newtheorem{corollary}{Corollary}

\def \Q {\ensuremath{\mathbb{Q}}}
\def \Aone {\ensuremath{A^{(1)}(p,p^{-s})}}
\def \Atwo {\ensuremath{A^{(2)}(p,p^{-s})}}
\def \Athree {\ensuremath{A^{(3)}(p,p^{-s})}}
\def \Afour {\ensuremath{A^{(4)}(p,p^{-s})}}

\def \Q {\mathbb{Q}}

\def \PG {\ensuremath{\mathcal{P}_G}}
\def \Wodpstratoff {\ensuremath{A_{\mbox{{\small ODP(strat)/off}}}(p,T)}}

\def \isom {\ensuremath{\stackrel{\sim}{=}}}

\def \T {\ensuremath{\mathfrak{T} }}
\def \Ttwo {\ensuremath{\mathfrak{T}_2 }}

\def \Qp {\mathbb{Q}_p}
\def \Z {\mathbb{Z}}

\def \Zp  {\mathbb{Z}_p}

\def \Fp  {\mathbb{F}_p}

\def \SF {\ensuremath{\mathcal{S}_F}}

\def \EFp {\ensuremath{|E(\mathbb{F}_p)|}}

\def \F23 {\ensuremath{F_{2,3}}}

\def \pthreetwo {\binom{3}{2}_p}
\def \ptwoone {\binom{2}{1}_p}

\begin{document}
\maketitle
\begin{abstract}
The local normal zeta functions of a finitely generated, torsion-free nilpotent group~$G$ of class $2$ depend on the
geometry of the Pfaffian hypersurface associated to the bilinear form induced by taking commutators in $G$.
 The smallest examples of zeta functions which are not finitely uniform arise from groups
  whose associated Pfaffian hypersurfaces are plane curves. In this paper we study groups
whose Pfaffians define singular curves, illustrating that the
local normal zeta functions may indeed invoke {\sl all} the degeneracy
loci of the Pfaffian.

\end{abstract}

\section{Introduction and statement of results}

In \cite{GSS/88} Grunewald, Segal and Smith defined the concept of the zeta function of a group. For a finitely
 generated, torsion-free nilpotent group $G$ -~a~\T-group, in short - they set
 $$a_n^\triangleleft:=|\{H\triangleleft G|\;|G:H|=n\}|.$$
 The group's {\sl normal zeta function} is the formal Dirichlet series
 $$ \zeta_{G}^\triangleleft(s):=\sum_{i=1}^\infty a_n^\triangleleft n^{-s}.$$
Here $s$ is a complex variable, and $\zeta_{G}^\triangleleft(s)$
defines an analytic function in some right
half-plane~$\{s\in\mathbb{C}|\,\mathfrak{Re}(s)>\alpha\}$ for some $\alpha\in\Q$. Similarly, the {\sl local} normal zeta function of $G$ at the prime~$p$ counts normal subgroups of $p$-power
index in $G$:
 $$\zeta_{G,p}^\triangleleft(s):=\sum_{n=0}^\infty a_{p^n}^\triangleleft p^{-ns}.$$
As $G$ is nilpotent, the local and global zeta functions are related by an {\sl Euler product}
$$\zeta_{G}^\triangleleft(s)=\prod_{p \text{ prime}}\zeta_{G,p}^\triangleleft(s) .$$
The local normal zeta functions are rational functions in $p^{-s}$. The global normal zeta function is
called {\sl finitely uniform} if there are finitely many rational functions $W_1(X,Y),\dots,W_n(X,Y)\in\mathbb{Q}(X,Y)$
 such that for every prime~$p$ there is an $r$ such that $\zeta^\triangleleft_{G,p}(s)=W_r(p,p^{-s})$, and {\sl uniform} if it is finitely
uniform with $n=1$.
(Finite) uniformity is not typical for zeta functions, however. The question of the dependence of the local factors on the
prime $p$ was linked to the classical problem of counting points on varieties mod~$p$ by du Sautoy and Grunewald. In \cite{duSG/00} they give a presentation of the local zeta functions as finite sums of products of rational functions
 in~$p$ and~$p^{-s}$ (which we will call the {\sl uniform components} of the
 presentation) and functions counting the number of $\Fp$-rational points of certain (boolean
 combinations of) algebraic varieties which are irreducible, smooth,
 and intersect normally. For an arbitrary collection of such
 varieties one would not expect these functions to be
 polynomials in~$p$. The description of varieties arising in this
 context has been an ongoing project ever since. We refer the reader
 to~\cite{duSSegal/00} for a survey on zeta functions of
 groups.

 The first example of a group whose zeta function is not
finitely uniform was given by du Sautoy~(\cite{duS-ecI/01},\cite{duS-ecII/01}). The local normal factors of this class-$2$-nilpotent group
 depend on the number of $\Fp$-points of the elliptic curve
$$ E=(y^2+x^3-Dx=0),\quad D\in\mathbb{N}.$$
He defined a nilpotent group
\begin{equation}
 G(E)=\langle x_1,\dots,x_6,y_1,y_2,y_3|\;[x_i,x_j]=M({\bf
   y})_{ij}\rangle,\label{group presentation}
\end{equation}
where $M({\bf y}):=\left(\begin{array}{cc}0&R({\bf y})\\-R({\bf
      y})^t&0\end{array}\right),$ and $$R({\bf y})=\left(\begin{array}{ccc}Dy_3&y_1&y_2\\y_1&y_3&0\\y_2&0&y_1\end{array}\right).$$
(Here we adopted additive notation for words in the $y_i$.) The curve~$E$
 appears as the {\sl Pfaffian hypersurface} 
$$\PG:=(\det(M({\bf y}))=0)$$ in~$\mathbb{P}^2(\Q)$ 
associated to this group. In~\cite{Voll/03} we proved

\begin{theorem} \label{elliptic-curve-theorem} For almost all primes~$p$
\begin{equation}
\zeta^\triangleleft_{G(E),p}(s)=W_1(p,p^{-s})+|E(\Fp)|W_2(p,p^{-s}),\label{e.c.
  theorem formula}
\end{equation}
where $|E(\Fp)|$ is the number of~$\Fp$-rational points of $E$ and
\begin{eqnarray*}
W_1(X,Y)&=&\frac{(1+X^{6}Y^{7}+X^{7}Y^{7}+X^{12}Y^{8}+X^{13}Y^{8}+X^{19}Y^{15})}{\prod_{i=0}^6(1-X^iY)\cdot(1-X^9Y^{18})(1-X^8Y^7)(1-X^{14}Y^8)},  \label{W1-ec} \\
W_2(X,Y)&=&\frac{(1-Y)(1+Y)X^6Y^5(1+X^{13}Y^8)}{\prod_{i=0}^6(1-X^iY)\cdot(1-X^9Y^{18})(1-X^8Y^7)(1-X^7Y^5)(1-X^{14}Y^8)}.  \label{W2-ec}
\end{eqnarray*}
\end{theorem}
Non-uniformity and the shape of the local zeta
functions, i.e. the existence of non-zero rational~$W_i(X,Y)$
in~(\ref{e.c. theorem formula}), were already established by du
Sautoy in~\cite{duS-ecII/01}. 

Moreover, we showed that in this manner (the number of $\Fp$-points of)
 every {\sl smooth} plane curve defined over~$\Q$ could be realised in local normal zeta functions of nilpotent
 groups. These improvements became feasible once we interpreted the
 local zeta function of a nilpotent group as a generating
   function associated to a weight function defined on the vertex set
 of the Bruhat-Tits building for $Sl_n(\Qp)$. In~\cite{Voll/03a} we
 generalised this approach to deal with \T-groups~$G$ of nilpotency
 class~$2$ - let us call them~$\Ttwo$-groups - whose Pfaffians~$\PG$ are
 smooth hypersurfaces containing no lines. This is the generic case
 for \Ttwo-groups~$G$ with `small centre and sufficiently large
 abelianization', or, more precisely, 
\begin{equation}
h(G')\leq 6 \text{  and  }h(G/{G'})>4
 h(G')-10.\label{lines condition}
\end{equation} 
(Here~$G'=[G,G]$ is the group's derived group and~$h(G)$ denotes the Hirsch length of the polycyclic
 group~$G$, i.e. the number of infinite cyclic factors in a
 decomposition series of~$G$.)

In this paper we treat classes of examples of \Ttwo-groups whose
Pfaffian~$\PG$ is singular. Our main Theorem~\ref{main theorem}
illustrates that in general the
formulae for the group's
local zeta functions may invoke not only (the number of $\Fp$-points
of) the Pfaffian hypersurface, but indeed all the strata of its rank stratification. 

For the rest of this paper $G$ shall denote a $\Ttwo$-group. Only for simplicity we assume
that $G'=Z(G)$, the group's centre, and $G/G'\cong \Z^d$, $d=2r$,
and $G'\cong \Z^{d'}$. Indeed, in general $Z(G)/G'$ is a finitely
  generated abelian group. As we are looking to prove results about
  almost all primes, we may assume it to be torsion-free. And as
 $$\zeta^\triangleleft_{G\times
      \Z^r}(s)=\zeta^\triangleleft_{G}(s)\cdot\zeta(s-n)\zeta(s-(n+1))\dots\zeta(s-(n+r-1))$$ where $h(G)=n$ and $\zeta(s)=\sum_{k=1}^\infty k^{-s}$ is the Riemann zeta function, we may assume $Z(G)=G'$. The condition ``$d$ even'' has also just been added for the convenience of presentation. We will assume that $G$ is presented as in~(\ref{group presentation}) above, where $M({\bf y})$ is an anti-symmetric $(d\times d)$-matrix of $\Z$-linear forms in
${\bf y}=(y_1,\dots,y_{d'})$ - but not necessarily in block form as in
the particular example~(\ref{group presentation}). We may then consider~$\PG$ as a hypersurface of
degree~$r$  in $\mathbb{P}^{d'-1}(\Q)$. We say that $G$ is {\sl
  regular} if
 $$\forall\,'\, p \mbox{ prime }\;\forall\,\alpha\in\mathbb{P}^{d'-1}(\Qp): \;rk(\overline{M(\alpha)})\geq 2(r-1).$$ 
(Here by $\forall\,'$ we mean `` for almost all'', i.e. for all but finitely
many.)

\begin{lemma}\label{lemma-regularity}
If $\PG\subseteq\mathbb{P}^{d'-1}(\Q)$ is smooth then~$G$ is regular.
\end{lemma}
\begin{proof} If $\PG$ is smooth, the reduction $\overline{\PG}$ is smooth for
  almost all primes~$p$. If $G$ were not regular there were infinitely
    many primes $p$ for which there exists a point ${\bf
    x}=(x_{ij})\in\mathbb{P}^{d'-1}(\Fp)$ for which $$rk(M({\bf
    x}))< 2(r-1).$$ It is thus enough to show that for almost all
    primes such a point must be a singular point of $\overline{\PG}$.

 Let~$\overline{i}$ be the reduction mod~$p$ of the linear embedding of~$\mathbb{P}^{d'-1}(\Q)$ into the
  projective space~$\mathbb{P}(\mathcal{S}_r)(\Q)$ of anti-symmetric
  $(2r\times2r)$-matrices over $\Q$ up to scalars given by
\begin{eqnarray*}
i:\mathbb{P}^{d'-1}(\Q)&\rightarrow& \mathbb{P}(\mathcal{S}_r)(\Q) \nonumber \\
\alpha& \mapsto &(\hdots :L_{ij}(\alpha):\dots) \nonumber
\end{eqnarray*}
where the linear forms $L_{ij}$ are given by the $ij$-th entry of the
matrix~$M(\alpha)$. Ignoring finitely many primes $p$ as we may
we assume that $\overline{i}$ is injective. Let $U:=\overline{i}(\mathbb{P}^{d'-1}(\Fp))$ and
$\mathcal{X}_r$ denote the universal Pfaffian hypersurface
in~$\mathbb{P}(\mathcal{S}_r)(\mathbb{Q})$ of {\sl singular}
anti-symmetric $(2r\times 2r)$-matrices up to scalars. Its singular locus $Sing(\mathcal{X}_r)$
consists of matrices of rank strictly less than~$2(r-1)$. We have~$\overline{\PG}=U\cap \overline{\mathcal{X}_r}$. It is enough to show that
$$ U\cap Sing(\overline{\mathcal{X}_r})\leq Sing(\overline{\PG}).$$
But that is clear since if ${\bf x}=(x_{ij})\in U\cap Sing(\overline{\mathcal{X}_r})$ then $T_{\bf x}(\overline{\PG})=U$. Thus~${\bf x}$ is a singular point of $\overline{\PG}$.
\end{proof}

The converse of Lemma~\ref{lemma-regularity} does not hold,
as the following example shows:

\begin{example}Consider the (singular) plane curve
$$ C:=(y^2-x^3+x^2=0). $$
We define two groups $G_1(C)$ and $G_2(C)$ as in~(\ref{group
  presentation}) above, with
$$
             R_1({\bf y})=\left( \begin{array}{ccc}
                                 y_1&y_2&0 \\
                                 0&y_1&y_2 \\
                                 -y_3&0&y_1-y_3
                                \end{array} \right),\;
R_2({\bf y})=\left( \begin{array}{ccc}
                                 y_1&y_2&0 \\
                                 0&y_1-y_3&y_2 \\
                                 -y_3&0&y_1
                                \end{array} \right).$$
Note that $C=\mathcal{P}_{G_1(C)}=\mathcal{P}_{G_2(C)}$, but that only
$G_2(C)$ is regular. As we will see, this distinction will be picked up by the
local normal zeta functions.
\end{example}

To state our main theorem, we have to introduce further notation and
 hypotheses. Assume for the rest of this chapter that $d'=3$, that
 $\sqrt{\det(M({\bf y}))}$ is non-zero and
 square-free, and that the Pfaffian curve~$\PG$ has irreducible components $E_i$,
$i=1,\dots,l$, all of degree~$>1$, and has at most ordinary double points (ODP's) as
 singularities, all of which are rational with rational slopes. Given a prime~$p$ we write
$$c_{p,i}:=|\overline{E}_i|$$ for the number of model points
over~$\Fp$ of~$E_i$, that is the number of $\Fp$-points of a
non-singular curve birationally equivalent to~$E_i$. For $j=1,2$ we write 
$$n_p^{(j)}:=|\{\alpha\in\mathbb{P}^2(\Fp)|\;\alpha \text{ is ODP of }
\overline{\PG}\text{ and }rk(M(\alpha))=2(r-j) \}|$$
for the number of $\Fp$-rational ODP's of {\sl rank deficit}~$2j$.
(In the above example we had $n_p^{(1)}=0$ and $n_p^{(2)}=1$ for
$G=G_1(C)$, $n_p^{(1)}=1$ and $n_p^{(2)}=0$ for
$G=G_2(C)$, and $c_{p,1}=1+p$ for both groups and all primes~$p$.)
With this notation we have

\begin{theorem}\label{main theorem} Assume that $G$ is a $\Ttwo$-group
  such that $\sqrt{\det(M({\bf y}))}$ is non-zero and square-free and the curve~$\PG$ has at most rational ordinary double
  points with rational slopes as
  singularities and has no lines. Then, for almost all primes~$p$,
\begin{eqnarray}
\frac{\zeta_{G,p}^\triangleleft(s)}{\zeta_{\Zp^d}(s)\zeta_{p}((d+3)s-3d)}&=&W_1(p,p^{-s})+\left(\sum_{i=1}^lc_{p,i}\right)W_2(p,p^{-s})\nonumber\\
&&\quad+n_p^{(1)}W_3(p,p^{-s})+n_p^{(2)}W_4(p,p^{-s}), \label{main
  theorem formula}
\end{eqnarray}
where 
\begin{eqnarray}
W_1(X,Y)&=&\frac{1+X^{2r}Y^{2r+1}+X^{2r+1}Y^{2r+1}+X^{4r}Y^{2r+2}+X^{4r+1}Y^{2r+2}+X^{6r+1}Y^{4r+3}}{(1-X^{4r+2}Y^{2r+2})(1-X^{2r+2}Y^{2r+1})},\nonumber\\&&\label{[1]} \\
W_2(X,Y)&=&\frac{(1-Y)(1+Y)X^{2r}Y^{2r-1}(1+X^{4r+1}Y^{2r+2})}{(1-X^{2r+1}Y^{2r-1})(1-X^{2r+2}Y^{2r+1})(1-X^{4r+2}Y^{2r+2})},\nonumber\\&& \label{[2]-[1]} \\
%the following is W_4 in thesis:
W_3(X,Y)&=&\frac{(1-Y)^2(1+Y)^2X^{2r}Y^{2r-3}(1+X^{2r+1}Y^{2r-1})(1+X^{4r+1}Y^{2r+2})}{(1-X^{2r+2}Y^{2r+1})(1-X^{2r}Y^{2r-3})(1-X^{2r+1}Y^{2r-1})(1-X^{4r+2}Y^{2r+2})},\nonumber\\&&\label{[1]+[4]-2[2]}
\\
%the following is W_5 in thesis:
W_4(X,Y)&=&\frac{(1-Y)(1+Y)X^{2r}Y^{2r-1}(1-X^{2r+2}Y^{2r-1})(1+X^{4r+1}Y^{2r+2})}{(1-X^{2r+1}Y^{2r-1})^2(1-X^{2r+2}Y^{2r+1})(1-X^{4r+2}Y^{2r+2})},\nonumber\\&&\label{[1]+[4a]-2[2]} 
\end{eqnarray}
\end{theorem}
In~(\ref{main theorem formula}), $\zeta_p(s)=(1-p^{-s})^{-1}$ denotes
the $p$-th local Riemann zeta function. The division on the left
hand side is purely formal.
\begin{corollary} \label{coro to main theorem}
Assume that, in addition, for each
$i=1,\dots,l$ there exist integers $B^{(i)}$ and complex numbers
$\beta_{i,j}$ such that
$$c_{p,i}=|\overline{E_i}|=1+p-\sum_{j=1}^{B^{(i)}}\beta_{i,j}$$ and a
$1-1$-correspondence
$$\left\{\frac{p}{\beta_{i,j}}| 1\leq j \leq
    B^{(i)}\right\}\stackrel{1-1}{\longleftrightarrow}\left\{\beta_{i,k}|
    1\leq k \leq  B^{(i)}\right\}.$$
Set 
\begin{eqnarray*}
c_{p^{-1},i}&:=&1+p^{-1}-\sum_{j=1}^{B^{(i)}}\beta_{i,j}^{-1}=p^{-1}c_{p,i}.\\
n_{p^{-1}}^{(j)}&:=&n_{p}^{(j)}.
\end{eqnarray*}
Then for all but finitely many primes $p$ the local normal zeta function~$\zeta^\triangleleft_{G,p}(s)$ satisfies a functional equation
\begin{equation}
\zeta^\triangleleft_{G,p}(s)|_{p \rightarrow p^{-1}}=-p^{\binom{2r+3}{2}-(4r+3)s}\zeta^\triangleleft_{G,p}(s).
\end{equation}
\end{corollary}

\begin{proofnodot} (of Corollary~\ref{coro to main theorem}): This follows immediately from the observation that
\begin{equation*}
W_i(p,p^{-s})|_{\substack{
                    p \rightarrow p^{-1} 
                 }}=\left\{\begin{array}{ll}
                                         p^{3}W_i(p,p^{-s})&\mbox{ for }i=1,3,4, \\
                                         p^{4}W_i(p,p^{-s})&\mbox{ for }i=2.
                                       \end{array} \right.
\end{equation*}
\end{proofnodot} %corollary to main theorem

Note that Theorem~\ref{elliptic-curve-theorem} follows from
Theorem~\ref{main theorem} by setting $r=3, l=1, c_1=\EFp,
n^{(1)}=n^{(2)}=0$, and that we get explicit formulae for the zeta
functions of the groups constructed in the example preceeding
Theorem~\ref{main theorem}. The conditions of
 Corollary~\ref{coro to main theorem} are satisfied in all these cases, giving us a local
functional equation.

Before we prove Theorem~\ref{main theorem} in Section~\ref{proof of main theorem}
 we give a short introduction to our method in Section~\ref{method} to make this paper self-contained. 

\begin{acknowledgements} This article was written when the author held
 a Postdoctoral Fellowship from the UK's Engineering and Physical
 Sciences Research Council (EPSRC). It comprises a part of the author's Cambridge PhD-Thesis, supported by the Studienstiftung des deutschen Volkes and the Cambridge European Trust. We would like to thank
 Marcus du Sautoy for his encouragement as supervisor, Fritz Grunewald
 for his continuous support and invaluable conversations, Burt Totaro
 for asking me about singular Pfaffians, and many
 others for sharing their knowledge and enthusiasm with me. 
\end{acknowledgements}

\section{Counting in Bruhat-Tits buildings} \label{method}
Now let $G$ be again a $\Ttwo$-group with $$Z(G)=G', h(G/G')=d=2r,
 h(G')=d'.$$
Nothing will be lost if we assume $d'\geq 2$, $d\geq4$.
For a fixed prime~$p$, the computation of the $p$-th normal local zeta function of~$G$
comes down to an emumeration of lattices (= subgroups of finite index)
in the $\Zp$-Lie algebra (with Lie brackets induced by taking commutators) $$G_p:=(G/Z(G)\oplus Z(G))\otimes_{\Z} \Zp$$ which
 are ideals in $G_p$. We
call a lattice $\Lambda\leq\Zp^n$ {\sl
  maximal} (in its homothety class) if $p^{-1}\Lambda\not\leq\Zp^n$. The key observation is the following

\begin{lemma}\label{gss-lemma}
For each lattice $\Lambda'\leq G_p'$ put $X(\Lambda')/\Lambda'=Z(G_p/\Lambda')$. Then
\begin{eqnarray*}
\lefteqn{\zeta^\triangleleft_{G,p}(s)=\zeta_{\Zp^{d}}(s)\sum_{\Lambda'\leq_f
  G_p'}|G_p':\Lambda'|^{d-s}|G_p:X(\Lambda')|^{-s}}\\
&=&\zeta_{\Zp^{d}}(s)\zeta_p((d+d')s-dd')\underbrace{\sum_{\substack{\Lambda'\leq_f G'_p  \\ \Lambda' {\text maximal}}}|G_p':\Lambda'|^{d-s}|G_p:X(\Lambda')|^{-s}}_{=:A(p,p^{-s}),{\text say}}
\end{eqnarray*}
\end{lemma}

This is essentially a local version of Lemma
6.1 in~\cite{GSS/88}. One of the main
theorems of~\cite{GSS/88} establishes the rationality of the
generating function\footnote{Recall
  that it is in general not true that $A(p,p^{-s})$ is a rational
  function in both $p^{-s}$ and $p$, i.e. that
  $\zeta^\triangleleft_{G}(s)$ is uniform.} $A(p,p^{-s})$
in~$p^{-s}$. Recall that (homothety classes of) maximal lattices are in one-to-one
correspondence with vertices of the Bruhat-Tits building $\Delta_{d'}$
for $Sl_{d'}(\Qp)$ (e.g.~\cite{Garrett/97}, \S 19).  To derive an
explicit formula for $A(p,p^{-s})$ requires a quantitative understanding of two weight functions defined on the set of vertices
of the building $\Delta_{d'}$. We write
\begin{equation}
A(p,p^{-s})=\sum_{[\Lambda']}p^{dw([\Lambda'])-sw'([\Lambda'])},\label{A}
\end{equation}
where, for a homothety class $[\Lambda']$ of a {\sl maximal} lattice
$\Lambda'$ in $G_p'\isom \Zp^{d'}$ we define
\begin{eqnarray*}
w([\Lambda']) &:=& \log_p(|G'_p:\Lambda'|),\\
w'([\Lambda']) &:=& w([\Lambda'])+ \log_p(|G_p:X(\Lambda')|)
\end{eqnarray*}

In \cite{Voll/03} we explain in detail how the evaluation of $w'$ might be
reduced to solving linear congruences. To recall the result in
Theorem~\ref{theorem linear congruences} below we have to introduce some
more notation. 

A maximal lattice $\Lambda'\leq\Zp^{d'}$ is said to be of {\sl type}
$\nu(\Lambda')=(I,{\bf r}_I)$ if $\Lambda'$ has elementary divisors
\begin{equation}
(p^{\nu}):=\left(\underbrace{1,\dots,1}_{i_1},\underbrace{p^{r_{i_1}},\dots,p^{r_{i_1}}}_{i_2-i_1},\dots,\underbrace
{p^{\sum_{j=1}^lr_{i_j}},\dots,p^{\sum_{j=1}^lr_{i_j}}}_{d'-i_l}\right)
\label{eldiv type}
\end{equation} for
\begin{equation}
I=\{i_1,\dots,i_l\}\leq\underline{d'-1}, 1\leq i_1<\dots<i_l\leq
 d'-1,\label{I}
\end{equation}
 and ${\bf r}_I=(r_{i_1},\dots,r_{i_l})\in \mathbb{N}_{>0}^l$.
Here we set~$\underline{n}:=\{1,\dots,n\}$. By abuse of notation we may say that a maximal lattice~$\Lambda'$ is
of type~$I$ if it is of type $(I,{\bf r}_I)$ for some positive vector~${\bf r}_I$, and that a (not necessarily maximal) lattice or a homothety
class has type~$I$ if the maximal element of the class has type $I$.

The group $\Gamma=Sl_n(\Zp)$ acts transitively on the set of maximal lattices of fixed type. If
we choose a basis for the $\Zp$-module $G_p'$ and represent lattices
as the row span of $d'\times d'$-matrices and denote by $\Gamma_{\nu}$
the stabilizer in $\Gamma$ of the lattice generated by the diagonal
matrix whose entries are given by the vector~(\ref{eldiv type}), the orbit-stabiliser
theorem gives us a one-to-one correspondence
\begin{equation}
\left\{\mbox{maximal lattices of type }(I,{\bf
    r}_I)\right\}\stackrel{1-1}{\longleftrightarrow}
    \Gamma/{\Gamma_{\nu}}.\label{correspondence}
\end{equation}

The correspondence~(\ref{correspondence}) allows us to describe 
$|G_p:X(\Lambda')|$ for a maximal lattice~$\Lambda'$ in terms of $M({\bf
  y})$, the matrix of commutators in a presentation for~$G$ as
in~(\ref{group presentation}). Its entries are $\Z$-linear forms in
generators $y_1,\dots,y_{d'}$ for $G'$.

\begin{theorem}\label{theorem linear congruences}[\cite{Voll/03}, \S 2.2] If~$\Lambda'$ corresponds
  to the coset~$\alpha\Gamma_{\nu}$ under~(\ref{correspondence}), where $\alpha\in\Gamma$ with
{\sl column vectors} $\alpha^j$, $j=1,\dots,d'$. Then the index
$|G_p:X(\Lambda')|$ equals the index of the kernel of the following system
of linear congruences in~$G_p/G_p'$:
\begin{equation}
\forall i \in \{1,\dots,d'\}\quad {\bf \overline{g}}M(\alpha^i)\equiv
0 \; mod \;(p^{\nu})_i,\label{linear congruences}
\end{equation}
where $ {\bf \overline{g}}=(\overline{g}_1,\dots,\overline{g}_d)\in
G_p/G'_p\isom \Zp^d$ and $(p^{\nu})_i$ denotes the $i$-th entry of the
vector $(p^{\nu})$ given in~(\ref{eldiv type}).
\end{theorem}

Further study of the linear congruences~(\ref{linear congruences}) will
enable us to formulate conditions on the group~$G$ under
which the explicit computation of $A(p,p^{-s})$ - and thus
$\zeta^\triangleleft_{G,p}(s)$ - is feasible. Given $\alpha\in\Gamma$ and~$I$ as in~(\ref{I}), let
$\overline{\alpha}$ denote the reduction mod~$p$ and define vector spaces
$$V_{i_j}:=\langle
\overline{\alpha^{i_j+1}},\dots,\overline{\alpha^{d'}}\rangle < \mathbb{P}^{d'}(\Fp)  ,\quad i_j\in I.$$
Clearly ${\text codim}(V_{i_j})=i_j-1$. We will call the nested sequence
of subspaces
$$ V_{i_1}>\dots>V_{i_l}$$
the {\sl flag of type~$I$ associated to $\alpha$} or indeed {\sl to the
maximal lattice $\Lambda'$} if $\nu(\Lambda')= I $ and $\Lambda'$
corresponds to $\alpha\Gamma_{\nu}$
under~(\ref{correspondence})\footnote{It is indeed straightforward to show that this
    is well-defined, i.e. independent of the coset representative
    $\alpha$.}. 
Given a fixed point $\overline{x}\in\mathbb{P}^{d'-1}(\Fp)$, we call a
lattice~$\Lambda'$ {\sl a lift of}~$\overline{x}$ if its associated
flag contains~$\overline{x}$ as $0$-dimensional member. Note that then
necessarily $d'-1\in\nu(\Lambda')$.

Given $I\leq\underline{d'-1}$ and $\overline{x}\in\mathbb{P}^{d'-1}(\Fp)$, we set

\begin{eqnarray}
A_I(p,p^{-s})&:=&\sum_{\nu(\Lambda')=I}p^{dw([\Lambda'])-sw'([\Lambda'])},\nonumber\\
A_{\overline{x}}(p,p^{-s})&:=&\sum_{\Lambda'\text{ lifts  }\overline{x}}p^{dw([\Lambda'])-sw'([\Lambda'])}.\nonumber
\end{eqnarray}

Clearly
\begin{eqnarray}
A(p,p^{-s})&=&\sum_{I\leq\underline{d'-1}}A_I(p,p^{-s}) \label{A rewrite 1}\\
&=&\sum_{\overline{x}\in\mathbb{P}^{d'-1}(\Fp)}A_{\overline{x}}(p,p^{-s})+\sum_{I\not\ni
  d'-1}A_I(p,p^{-s})\label{A rewrite 2}
\end{eqnarray}

Formula~(\ref{A rewrite 2}) gives us an inkling how the combinatorial
object $A(p,p^{-s})$ is linked to the geometry of the Pfaffian~$\PG$,
a hypersurface in projective $(d'-1)$-space. Theorem~\ref{theorem linear congruences} tells us that the number
$|G_p:X(\Lambda')|$ depends in the first instance on the way the flag
associated to $\Lambda'$ intersects with~$\overline{\PG}$. Assume from
now on that its defining polynomial~$\sqrt{\det(M({\bf y}))}$ is
non-zero and square-free, and that $\PG$ contains no lines. This
last condition will be extremely useful and is generically
satisfied if $r>2d'-5$ (see Beauville's Appendix to~\cite{Voll/03a} for details). 

In the easiest case $\PG$ has no $\Fp$-rational points. Thus
$\det(M(\alpha))$ is a $p$-adic unit for all
$\alpha\in\Zp^{d'}\setminus(p\Zp^{d'})$ and~(\ref{linear congruences}) is
equivalent to the single congruence
$$ {\bf \overline{g}} \equiv 0 \; mod \; p^{\sum_{j=1}^l
  r_{i_j}}.$$ 
Hence 
\begin{equation}
|G_p:X(\Lambda')|=p^{d\sum_{j=1}^l
  r_{i_j}}. \label{index}
\end{equation}
In particular, both weight functions $w$ and $w'$ only
depend on the maximal lattices' {\sl type}.
Such functions have been studied in~\cite{Voll/03a}:

\begin{theorem}\label{theorem igusa}[\cite{Voll/03a}] Assume that the weight functions $w([\Lambda'])$ and $w'([\Lambda'])$ only depend on the
elementary divisor type of~$\Lambda'$. Then there exist integers
$a_i,b_i$ for $i\in\{1,\dots,d'-1\}$, depending only on $d,d'$ and
$i$, such that, for all $I\leq \underline{d'-1}$
\begin{equation}
A_I(p,p^{-s})=b_I(p^{-1})\prod_{i\in
  I}\left.\frac{X_i}{1-X_i}\right|_{X_i\rightarrow
  p^{b_i-a_is}}\label{igusa's function}
\end{equation}
where $b_I(p)\in\Z[p]$ denotes the number of $\Fp$-rational points of
the projective variety of flags of type $I$ (and
$b_{\emptyset}(p)=1$). 
In particular
\begin{equation}
  A_{\overline{x}}(p,p^{-s})=\left.\underbrace{\frac{p^{-(d'-1)}X_{d'-1}}{1-X_{d'-1}}}_{=:\Aone,\text{
  say}}\left(\sum_{I'\in\underline{d'-2}}b_{I'}(p^{-1})\prod_{i\in
  I'}\frac{X_i}{1-X_i}\right)\right|_{X_i\rightarrow
  p^{b_i-a_is}}, \label{A_x non-rational point}
\end{equation}
independently of $\overline{x}\in\mathbb{P}^{d'-1}(\Fp)$. The following local functional equation holds
$$A(p,p^{-s})|_{p\rightarrow p^{-1}}=(-1)^{d'-1}p^{-\binom{d'}{2}}A(p,p^{-s}).$$
\end{theorem}
In this special case $A(p,p^{-s})$ is an instantiation of a type of
functions first studied by Igusa~\cite{Igusa/89}. The polynomials
$b_I(p)$ are easily computed in terms of $q$-binomial polynomials (cf.~\cite{Macdonald/95}). 
The proof of the functional equation relies on a symmetry of the
(poset) lattice of Schubert cells for the complete flag variety. We
note that this symmetry is independent of the `numerical data'~$a_i,b_i$. See~\cite{Voll/03a}, \S 4.2.1, for details. There it is
also proved that 
\begin{equation}
a_i:=d+d'-i,\; b_i:=(d+i)(d'-i), \;1\leq i \leq d'-1\label{numerical data}
\end{equation}
 is the right numerical data for
the case that $\PG$ has no $\Fp$-points. Formula~(\ref{[1]}) may then be
retrieved from~(\ref{A rewrite 1}), (\ref{igusa's function}) and~(\ref{numerical data}).

Note that under the assumption that $\PG$ contains no
higher-dimensional linear subspaces, the formulae for the
$A_I(p,p^{-s})$, $d'-1\not\in I$, are given by~(\ref{igusa's
  function}) regardless if $\PG$ has $\Fp$-rational points, at least for almost all primes~$p$. This follows from
inspection of the linear congruences~(\ref{linear congruences}). It can also be shown
(\cite{Voll/03a}, Prop. 3) that under this assumption, for
$\overline{x}\in\mathbb{P}^{d'-1}(\Fp)$, for almost all primes,
\begin{equation}
A_{\overline{x}}(p,p^{-s})=\underbrace{\left(\sum_{\substack{\Lambda'\text{ lifts }\overline{x}\\\nu(\Lambda')=d'-1}}p^{dw([\Lambda'])-sw'([\Lambda'])}\right)}_{=:A'_{\overline{x}}(p,p^{-s})}\left.\left(\sum_{I'\in\underline{d'-2}}b_{I'}(p^{-1})\prod_{i\in
  I'}\frac{X_i}{1-X_i}\right)\right|_{X_i\rightarrow
  p^{b_i-a_is}}\label{A_x rational point}
\end{equation}

If $\overline{x}\not\in\overline{\PG}$ and the numerical data is as
in~(\ref{numerical data}) above, this formula coincides with
the one given in~(\ref{A_x non-rational point}). If $\overline{x}$ is
a rational point of $\overline{\PG}$, however, equation~(\ref{A_x rational point}) tells
us that the mistake we have made by choosing~(\ref{A_x non-rational
  point}) instead is localised in the first factor. To correct it we
have to compute $A'_{\overline{x}}(p,p^{-s})$, a sum which extends just over
lattices lifting the (highly incomplete) flag consisting just of the
point $\overline{x}$. This generating function will depend on the structure
of $\PG(\Qp)$ near $\overline{x}$. The easiest case - that
$\overline{x}$ is a smooth point of $\overline{\PG}$ - was dealt with
in~\cite{Voll/03a} and will be recalled in the remainder of this
chapter. The main point of Theorem~\ref{main theorem} is to illustrate
that to deal with even the tamest singularities requires additional
information (here on the rank filtration) and intricate
computations. A satisfactory machinery to compute
$A'_{\overline{x}}(p,p^{-s})$ in general would be highly desirable but
seems out of reach by far. It would have to implement - among
other things - a resolution of singularities for~$\PG$.
 We will study the effects of dropping the simplifying assumption ``$\PG$ contains no lines'' in a future paper.

If $\overline{x}$ is a rational point of~$\overline{\PG}$, the linear
congruences~(\ref{linear congruences}) to be solved in order to
compute~$A'_{\overline{x}}(p,p^{-s})$ are
\begin{eqnarray}
 {\bf \overline{g}}M(\alpha^{d'})&\equiv& 0 \; \mod \; p^{r_{d'-1}}\label{congruences2},
\end{eqnarray}
where $\alpha^{d'}\in\mathbb{P}^{d'-1}(\Z/(p^{r_{i_l}}))$ is a lift of
$\overline{x}$. If $\overline{x}$ is a smooth point of $\PG$ we can find local
coordinates such that~(\ref{congruences2}) reads like

\begin{eqnarray}
{\bf \overline{g}}\mbox{ diag}\left( \left(\begin{array}{cc}0 & x \\
                                                             -x&0
                                                           \end{array}\right),J_2,\dots,J_2\right) & \equiv & 0 \; \mod p^{r_{d'-1}},\label{congruences2A}
\end{eqnarray}
where $J_2=\left(\begin{array}{rl} 0&1\\-1&0\end{array} \right)$ and
$x\in p\Zp/(p^{d'-1})$. Thus
\begin{equation}
|G_p:X(\Lambda')|=p^{d\sum_{j=1}^l r_{i_j}-2\min\{r_{d'-1},v_p(x)\}} \label{index-smooth},
\end{equation}
with $p$-adic valuation $v_p$.

A quantitative version of Hensel's
Lemma describes how $v_p(x)$ varies as~$\Lambda'$ runs over the maximal lattices lifting $\overline{x}$ and
allows us to compute the generating function
$$A'_{\overline{x}}(p,p^{-s})=\frac{p^{2r-s(2r-1)}(1-p^{2r+1-s(2r+1)})}{(1-p^{2r+1-s(2r-1)})(1-p^{2r+2-s(2r+1)})}=\Atwo,\text{
  say,}$$ and thus the ``correction term''

\begin{equation}
W_2(p,p^{-s})=\left.\left(\Atwo-\Aone\right)\left(\sum_{I'\in\underline{d'-2}}b_{I'}(p^{-1})\prod_{i\in
  I'}\frac{X_i}{1-X_i}\right)\right|_{X_i\rightarrow
  p^{b_i-a_is}},\nonumber
\end{equation}
giving us~(\ref{[2]-[1]}) explicitly.  See~\cite{Voll/03a}, \S 4.2, for details. A variant of this Lemma will be derived in
Section~\ref{proof of main theorem} (formula~(\ref{preimage phi -
  smpt})) to deal with singular points.

%\sum_{\Lambda'\text{ lifts
%  }\overline{x}}p^{w([\Lambda'])}\left(p^{-s
%    w'([\Lambda'])}-p^{-s\left(w([\Lambda'])+d\sum_{j=1}^{l}r_{i_j}\right)} \right).
%.......................................................................
%....................................................................
\section{Proof of Theorem~\ref{main theorem}}\label{proof of main theorem}

To prove Theorem~\ref{main theorem} we have to compute the two
`correction terms' $W_3(p,p^{-s})$ and $W_4(p,p^{-s})$. We have seen
that it suffices to compute $A'_{\overline{x}}(p,p^{-s})$ for a fixed
point of rank deficit~$2$ and $4$, respectively. Note that $d'=3$.

%...............................................................................
\subsection{$W_4(p,p^{-s})$}
%.............................................................................

% ODP strat, in Thesis-notation
Let $\overline{x}$ be a fixed ODP of~$\overline{\PG}$ of rank
deficit~$4$. Dismissing finitely many primes if necessary, we may
assume that we can find local affine coordinates such that condition~(\ref{congruences2}) reads like

\begin{equation}
{\bf \overline{g}} \mbox{ diag}\left( \left(\begin{array}{cc}0 & x_1 \\
                                                             -x_1&0
                                                           \end{array}\right),\left(\begin{array}{cc}0 & x_2 \\
                                                             -x_2&0
                                                           \end{array}\right),J_2\right) \equiv 0 \; \mod p^{a},\label{odpstrat eqns 2}
\end{equation}
($x_1,x_2\in p\Z/(p^{a})$, $a=r_2\geq 1$). We compute the weight function as
\begin{equation}
\boxed{w'([\Lambda'])={(d+1)a-2\left(\min\{a,v_p(x_1)\}+\min\{a,v_p(x_2)\}\right)}.} \nonumber
%\label{weight function:A_odp(strat)/off}
\end{equation}
We
shall factorize the map
 $$[\Lambda'] \mapsto p^{dw([\Lambda'])-sw'([\Lambda'])}$$ 
 over the set $N:=\{(a,b,c) \in \mathbb{N}_{>0}^3|\;a\geq b, a \geq
 c\}$ as $\psi \circ \phi$ where 
\begin{eqnarray}
\phi:[\Lambda'] &\mapsto& (a,v_p(x_1),v_p(x_2))\nonumber \\
\psi:(a,b,c) &\mapsto& p^{2ra-s((2r+1)a-2\min(a,b)-2\min(a,c))} \label{psi}
\end{eqnarray}
To deduce an explicit formula for the rational generating
function~$A_{\overline{x}}(p,p^{-s})$ we must study the
cardinalities of the fibres of $\phi$ and must brake up the cone~$N$ in
order to eliminate the ``$\min$'' in~(\ref{psi}). It is indeed not hard to compute
\begin{equation} \label{preimage phi - smpt} 
|\phi^{-1}(a,b,c)|=\left\{\begin{array}{ll}
                  1 & \mbox{if }a=b=c, \\
                  (1-p^{-1})\;p^{a-b} & \mbox{if }a>b, a=c, \\
                  (1-p^{-1})\;p^{a-c} & \mbox{if }a>c, a=b, \\
                  (1-p^{-1})^2p^{2a-b-c} & \mbox{if }a>b, a>c,
                \end{array} 
                \right.
\end{equation}
and to see that the following decomposition of~$N$ suits our needs:
\begin{eqnarray}
N&=&N_0+N_1+N_2+N_3,\text{ where }\label{N into four parts} \\
N_0&:=&\{(a,b,c)\in N|\;a=b=c\},\nonumber\\
N_1&:=&\{(a,b,c)\in N|\;a=c>b\},\nonumber \\
N_2&:=&\{(a,b,c)\in N|\;a=b>c\},\nonumber \\
N_3&:=&\{(a,b,c)\in N|\;a>b,\,a>c\}. \nonumber 
\end{eqnarray}
Indeed, if we set $n_j:=\dim(N_j)$, we have 
\begin{eqnarray}
\lefteqn{(a,b,c)\in N_j \Rightarrow}\nonumber\\
&&\frac{|\phi^{-1}(a,b,c)|\psi(a,b,c)}{(1-p^{-1})^{n_j-1}}=m_{j,X}(p,p^{-s})^am_{j,Y}(p,p^{-s})^bm_{j,Z}(p,p^{-s})^c \label{Laurent}
\end{eqnarray}
for suitable Laurent monomials\footnote{ A {\sl Laurent monomial} in
  variables $x_1,\dots,x_n$ is a
  power product $x_1^{r_1}\dots x_n^{r_n}$ with $r_i\in\Z$.} $m_{j,X}(p,p^{-s}),m_{j,Y}(p,p^{-s})$
and $m_{j,Z}(p,p^{-s})$ in variables~$p$ and $p^{-s}$.

 Let $F_j(X,Y,Z)$, $j=0,\dots,3$ be the zeta function of the cone~$N_j$, i.e. the rational function 
$$F_j(X,Y,Z)=\sum_{\substack{(a,b,c)\in N_j}}X^aY^bZ^c.$$

\noindent Simply substituting and summing up yields
\begin{eqnarray}
\lefteqn{\Afour:=A'_{\overline{x}}(p,p^{-s})=\sum_{\substack{\Lambda'\text{ lifts
    }\overline{x}\\\nu(\Lambda')=2}}p^{dw([\Lambda'])-sw'([\Lambda'])}=}\nonumber \\
&&\sum_{j=0}^3\sum_{(a,b,c)\in N_j}|\phi^{-1}(a,b,c)|\psi(a,b,c)=  \left.\sum_{j=0}^3(1-p^{-1})^{n_j-1}F_j(X,Y,Z)\right|_{\substack{X=m_{j,X}(p,p^{-s})\\Y=m_{j,Y}(p,p^{-s})\\Z=m_{j,Z}(p,p^{-s})}}=\nonumber\\
&&\frac{p^{2r-s(2r-3)}(1+p^{2r+1-s(2r-1)}-2p^{2r-s(2r-1)}-2p^{2r-1-s(2r+1)}+p^{4r+1-4sr}+p^{2r-s(2r+1)})}{(1-p^{2r-s(2r-3)})(1-p^{2r+1-s(2r-1)})(1-p^{2r+2-s(2r+1)})},\nonumber
\end{eqnarray}

\noindent Table~\ref{odp(strat)-table} records the generating functions $F_j(X,Y,Z)$
together with the integers~$n_j$ and the Laurent
monomials~$m_{j,X}(p,p^{-s})$, $m_{j,Y}(p,p^{-s})$,
$m_{j,Z}(p,p^{-s})$. The latter are easily read off from~(\ref{psi})
and~(\ref{preimage phi - smpt}). Routine computations with rational
functions now yield the desired expression~(\ref{[1]+[4a]-2[2]}) for

\begin{eqnarray}
\lefteqn{W_4(p,p^{-s})=}\\
&&\left( \Afour - 2\Atwo +
  \Aone\right)\left.\left(1+(p^{-1}+1)\frac{X_1}{1-X_1}\right)\right|_{X_1\rightarrow p^{4r+2-s(2r+2)}} \nonumber
\end{eqnarray}

\begin{table}[h]
\centering
\setlength{\extrarowheight}{7pt}
\begin{tabular}{c|c|c|c|c|c}
j&$n_j$&$F_j(X,Y,Z)$&$m_{j,X}(p,p^{-s})$&$m_{j,Y}(p,p^{-s})$&$m_{j,Z}(p,p^{-s})$ \\ \hline 
1&1&$\frac{XYZ}{1-XYZ}$&$p^{2r-s(2r+1)}$&$p^{-2s}$&$p^{-2s}$ \\
2&2&$\frac{X^2YZ^2}{(1-XYZ)(1-XZ)}$&$p^{2r+1-s(2r+1)}$&$p^{-1-2s}$&$p^{-2s}$  \\
3&2&$\frac{X^2Y^2Z}{(1-XYZ)(1-XY)}$&$p^{2r+1-s(2r+1)}$&$p^{-2s}$&$p^{-1-2s}$ \\
4&3&$\frac{X^2YZ(1-X^2YZ)}{(1-XYZ)(1-XY)(1-XZ)(1-X)}$&$p^{2r+2-s(2r+1)}$&$p^{-1-2s}$&$p^{-1-2s}$ 
\end{tabular}
\caption{Computing $\Afour$.}
\label{odp(strat)-table}
\end{table}

%............................................................................
\subsection{$W_3(p,p^{-s})$}
%.......................................................................

Now let $\overline{x}$ be a fixed ODP of~$\overline{\PG}$ of rank deficit~$2$. Ignoring
finitely many primes as we may, we choose
local coordinates such that condition~(\ref{congruences2}) becomes
\begin{equation}
{\bf \overline{g}} \mbox{ diag}\left( \left(\begin{array}{cc}0 & x_1x_2 \\
                                                             -x_1x_2&0
                                                           \end{array}\right),J_2,J_2\right) \equiv 0 \; \mod p^{a},\label{odpnonstrat eqns 2}
\end{equation}
($x_1,x_2\in p\Z/(p^{a})$, $a=r_2\geq1$). We compute the weight function
\begin{equation*}
\boxed{w'([\Lambda'])={(d+1)a-2\min\{a,v_p(x_1)+v_p(x_2)\}}.} \label{w' for flagtype ODP(nonstrat) - offline}
\end{equation*}

\noindent Again we shall  factorize
 $$[\Lambda'] \mapsto p^{dw([\Lambda'])-sw'([\Lambda'])}$$ 
 over the set $N:=\{(a,b,c) \in \mathbb{N}_{>0}^3|\;a\geq b, a \geq c\}$ as $\psi\circ \phi$ where 
\begin{eqnarray*}
\phi:[\Lambda'] &\mapsto& (s,v_p(x_1),v_p(x_2)) \\
\psi:(a,b,c) &\mapsto& p^{2ra-s(2r+1)a-2\min(a,b+c)}
\end{eqnarray*}
Recall the formula for $|\phi^{-1}(a,b,c)|$ given in~(\ref{preimage
  phi - smpt}). The subdivision~(\ref{N into four parts}) of $N$,
  however, has to be refined to
\begin{eqnarray}
N&=&\sum_{i=0}^5N_i,,\text{ where } N_0,N_1,N_2 \text{ as in (\ref{N
    into four parts}) and } \label{N into six parts} \\
N_3&:=&\{(a,b,c)\in N|\;a=b+c\}, \nonumber \\
N_4&:=&\{(a,b,c)\in N|\;a>b+c\}, \nonumber \\
N_5&:=&\{(a,b,c)\in N|\;a<b+c\}. \nonumber 
\end{eqnarray}

If we set $n_j:=dim(N_j)$ for $j=0,1,2$ and $n_j:=3$ for $j=3,4,5$,
equation~(\ref{Laurent}) holds again all $j$ and suitable Laurent monomials,
recorded in~Table~\ref{odp(nonstrat)-table}.

\begin{table}[h]
\centering
\setlength{\extrarowheight}{9pt}
\begin{tabular}{c|c|c|c|c|c}
j&$n_j$&$F_j(X,Y,Z)$&$m_{j,X}(p,p^{-s})$&$m_{j,Y}(p,p^{-s})$&$m_{j,Z}(p,p^{-s})$ \\ \hline 
0&1&$\frac{XYZ}{1-XYZ}$&$p^{2r-s(2r-1)}$&$1$&$1$ \\
1&2&$\frac{X^2YZ^2}{(1-XYZ)(1-XZ)}$&$p^{2r+1-s(2r-1)}$&$p^{-1}$&$1$ \\
2&2&$\frac{X^2Y^2Z}{(1-XYZ)(1-XY)}$&$p^{2r+1-s(2r-1)}$&$1$&$p^{-1}$ \\
3&3&$\frac{X^2YZ}{(1-XY)(1-XZ)}$&$p^{2r+2-s(2r-1)}$&$p^{-1}$&$p^{-1}$ \\
4&3&$\frac{X^3YZ}{(1-XY)(1-XZ)(1-X)}$&$p^{2r+2-s(2r+1)}$&$p^{-1-2s}$&$p^{-1-2s}$ \\
5&3&$\frac{X^3Y^2Z^2}{(1-XY)(1-XZ)(1-XYZ)}$&$p^{2r+2-s(2r-1)}$&$p^{-1}$&$p^{-1}$ \\
\end{tabular}
\caption{Computing $\Athree$.}
\label{odp(nonstrat)-table}
\end{table}
 
Substituting and summing up yields
\begin{eqnarray}
\lefteqn{\Athree:=A'_{\overline{x}}(p,p^{-s})=\sum_{\substack{\Lambda'\text{ lifts
    }\overline{x}\\\nu(\Lambda')=2}}p^{dw([\Lambda'])-sw'([\Lambda'])}=}\nonumber\\
&&\sum_{j=0}^{5}\sum_{(a,b,c)\in N_j}|\phi^{-1}(a,b,c)|\psi(a,b,c)=\left.\sum_{j=0}^{5}(1-p^{-1})^{n_j-1}F_j(X,Y,Z)\right|_{\substack{X=m_{jX}(p,p^{-s})\\Y=m_{jY}(p,p^{-s})\\Z=m_{jZ}(p,p^{-s})}}
=\nonumber\\
&&\frac{p^{2r-s(2r-1)}(1-2p^{2r+1-s(2r-1)}+p^{2r+2-s(2r-1)}-p^{2r+2-s(2r+1)}+p^{4r+2-s4r})}{(1-p^{2r+1-s(2r-1)})^2(1-p^{2r+2-s(2r+1)})} 
\nonumber\\&& \label{gen.fu.odp(nonstrat)}
\end{eqnarray}
and finally the desired formula~(\ref{[1]+[4]-2[2]}) for
\begin{eqnarray}
\lefteqn{W_3(p,p^{-s})=}\\
&&\left( \Athree - 2\Atwo +
  \Aone\right)\left.\left(1+(p^{-1}+1)\frac{X_1}{1-X_1}\right)\right|_{X_1\rightarrow p^{4r+2-s(2r+2)}}. \nonumber
\end{eqnarray}
See Table~\ref{odp(nonstrat)-table} for the zeta functions of the
respective cones and the corresponding Laurent monomials.

\bibliographystyle{amsplain}
\bibliography{thebibliography}

\end{document}